\documentclass[12pt]{article}
\usepackage{color}   
\usepackage{amsmath, amssymb}
\usepackage[utf8]{inputenc}
\usepackage[dvips]{graphicx}
\usepackage{units} 
\usepackage{t1enc}
\usepackage{bbm}
\newtheorem{lem}{Lemma}
\newtheorem{dfn}{Definition}

\newtheorem{prop}{Proposition}
\newtheorem{thm}{Theorem}

\newcommand{\bee}{\begin{equation}}
\newcommand{\eee}{\end{equation}}
\newcommand{\bees}{\begin{equation*}}
\newcommand{\eees}{\end{equation*}}
\newcommand{\bali}{\begin{aligned}}
\newcommand{\eali}{\end{aligned}}
\newcommand{\dy}{\:\mathrm{d}y\:}

\newcommand{\dx}{\:\mathrm{d}x\:}
\newcommand{\ds}{\:\mathrm{d}s\:}
\newcommand{\pa}{\partial}
\newcommand{\er}{\mathbb{R}}
\newcommand{\erp}{\mathbb{R}^+}

\newcommand{\mcd}{\mathcal{D}}
\newcommand{\mcn}{\mathcal{N}}
\newcommand{\mcf}{\mathcal{F}}

\newcommand{\uu}{\mathbf{u}}
\newcommand{\ee}{\mathbf{e}}
\newcommand{\vv}{\mathbf{v}}

\title{Numerical solution of fractional order diffusion problems 
with Neumann boundary conditions}

\author{B\'ela J. Szekeres$^{1}$
and Ferenc Izs\'ak$^{1,}$
\footnote{Corresponding author. E-mail address:
\texttt{izsakf@cs.elte.hu},
tel.: +36 13722500-8428, fax: +36 13812158.}}
   
\begin{document}

\maketitle
\begin{center}
       $^1$ Department of Applied Analysis and Computational Mathematics, \\
       E\"otv\"os Lor\'and University
       H-1117, Budapest, P{\'a}zm{\'a}ny P. s. 1/C, Hungary.
\smallskip
\end{center}

\begin{abstract}
A finite difference numerical method is investigated for fractional 
order diffusion problems in one space dimension.
For this, a mathematical model is developed to incorporate homogeneous 
Dirichlet and Neumann type boundary conditions. 
The models are based on an appropriate extension
of the initial values. The well-posedness of the obtained initial
value problems is proved and it is pointed out that the extensions 
are compatible with the above boundary conditions.
Accordingly, a finite difference scheme is constructed for the Neumann
problem using the shifted Gr\"unwald--Letnikov approximation of the 
fractional order derivatives, which is based on infinite many basis points.  
The corresponding matrix is expressed in a closed form and the 
convergence of an appropriate implicit Euler scheme is proved.
\end{abstract}

\section{Introduction}
A widely accepted constitutive relation, the first Fick's law leads to 
the standard diffusion model. At the same time, more observations 
confirmed the presence of super- and subdiffusive dynamics in several 
phenomena. They range from the plasma physics \cite{treumann97}
through population dynamics \cite{edwards07} and groundwater flows 
\cite{benson00} to the anomalous diffusion of some chemical compounds. 
These observations inspired the application of the 
fractional calculus, which has a long history \cite{podlubny99}. 
An alternative approach for modeling superdiffusion can be given in the 
framework of the nonlocal calculus, which has been recently developed, 
see \cite{du13} and \cite{du12} for an up-to-date overview with further
references.

Many attempts were made for the numerical solution of the corresponding 
space-fractional PDE's. Most of them based on finite difference
discretization. It was pointed out that a non-trivial discretization
of the one sided fractional derivatives lead to a stable method 
\cite{meerschaert04}. This result has been generalized in many 
aspects: higher-order methods and multi-dimensional schemes
\cite{deng14},\cite{tadjeran07},\cite{tadjeran06},\cite{tian12},\cite{zhou13} 
were constructed and analyzed. 
Note that recently the finite element (Galerkin) discretization was
initiated \cite{nochetto14} for the fractional diffusion equations
and a composite approach was analyzed \cite{huang13}.

The problem which inspired the present research is the following. 
If a superdiffusive evolution of some density $u$ is observed on a 
physical volume then we have information on the density only on
the closure of this. On the other hand, in the mathematical model, 
nonlocal operators are used, which require the value of $u$
also outside of the domain. For an accurate finite difference 
approximation we also need values outside of the domain. 
In this way, it is natural to look for an extension of the density $u$. 

The majority of the authors consider homogeneous Dirichlet boundary 
conditions and assume \emph{zero} values outside of the domain. 
A similar approach is applied in the non-local calculus \cite{du12}. 

Regarding Neumann type boundary conditions, at our best knowledge, there
is only one attempt \cite{ilic05}, dealing the boundary condition at the 
operator level and proposing a matrix transformation technique.    

The aim of this paper is to 
\begin{itemize}
\item
develop a meaningful mathematical approach to model homogeneous Neumann 
(and Dirichlet) boundary conditions  
\item
analyze the well-posedness of the corresponding PDE's
\item
construct a corresponding finite difference scheme with a full error analysis.
\end{itemize}

\section{Mathematical preliminaries}
\paragraph{Fractional calculus}
We summarize some basic notions and properties of the fractional 
calculus. For more details and examples we refer to the monographs
\cite{kilbas06},\cite{miller93},\cite{samko93}.\\
To define an appropriate function space for the fractional order 
derivatives on the real axis we introduce for arbitrary $a,b\in\er$
the function spaces 
$$
\bar C(a,b) = C[a,b]|_{(a,b)}
\quad\textrm{and}\quad
\bar C(a,b)/\er = \{ f\in \bar C(a,b): \int_a^b f = 0\}. 
$$
With these we define
$$
C_I(\er) = \{\textrm{$b-a$ - periodic extension of $f$: 
$f \in\bar C(a,b)/\er\}$}. 
$$
\emph{Remarks:}\\
1. Functions in $C_I(\er)$ are bounded and 
they are continuous except of the possible discontinuity points 
$\{a + k(b-a):\;k\in\mathbb{Z}\}$.\\
2. In the points $\{a + k(b-a):\;k\in\mathbb{Z}\}$ we define 
$f$ to be $\frac{f(a)+f(b)}{2}$.  
\quad$\square$

\begin{dfn}\label{fr_der}
For the exponent $\beta\in (0,1)$ the fractional order integral 
operators $_{-\infty}^{}I_x^\beta$ and $_{x}^{}I_\infty^\beta$ on the space 
$C_I(\er)$ are defined with 
$$
_{-\infty}^{}I_x^\beta f(x) = 
\frac{1}{\Gamma(\beta)}\int_{-\infty}^{x} \frac{f(s)}{(x-s)^{1-\beta}}\:\mathrm{d}s
$$
and 
$$
_{x}^{}I_\infty^\beta f(x) = 
\frac{1}{\Gamma(\beta)}\int_{x}^{\infty} \frac{f(s)}{(s-x)^{1-\beta}}\:\mathrm{d}s.
$$
With this, the left and right-sided \emph{Riemann--Liouville derivatives} 
of order $\alpha\in\erp\setminus\mathbb{Z}$ are given by
$$
_{-\infty}^{RL} \pa_x^{\alpha} f(x) =
\pa^n_x{}\: _{-\infty}^{}I_x^{n-\alpha} f(x)
$$
and 
$$
_x^{RL}\pa_\infty^\alpha f(x) =
(-1)^n \pa^n_x{} \:_x^{}I_\infty^{n-\alpha} f(x),
$$
where $n$ is the integer with $n-1<\alpha < n$.\\
Accordingly, for a function $f = f_0 + C_f$, where $f_0 \in C_I(\er)$ 
and $C_f$ is a constant function we define the Riesz derivative with
$$
\pa_{|x|}^\alpha f(x) = 
  C_{\sigma} (_{-\infty}^{RL}\pa_x^\alpha f_0(x) 
+ _x^{RL}\pa_\infty^\alpha f_0(x)),
$$
where $C_{\sigma} = -\frac{\sigma}{2\cos\alpha\frac{\pi}{2}}$ with 
a given positive constant $\sigma$.
\end{dfn} 
\noindent\emph{Remarks:}\\
1.  For the simplicity, the constant $\sigma$ -- which corresponds 
to the intensity of the superdiffusive process in the real-life 
phenomena -- does not appear in the notation $\pa_{|x|}^\alpha$.\\
2. In the original definition, the Riemann--Liouville derivatives 
are given on a bounded interval $(a,b)\subset\er$ such that in the above 
definition $-\infty$ and $\infty$ are substituted with $a$ and $b$, 
respectively. 
In this case $_{a}^{}I_x^\beta$ and $_{x}^{}I_b^\beta$ can be defined
for the exponent $\beta\in\erp$ or even for $\beta\in\mathbb{C}$
with $\textrm{Re}\:\beta>0$, in case of complex valued functions, see 
Section 2.1 in \cite{kilbas06}.
Moreover, the definition can be extended to be a bounded operator
on the function space $L_{1}(a,b)$, see Theorem 2.6 in \cite{samko93}.
For an overview of the alternating notations and definitions, we 
refer to the review paper \cite{hilfer08}.\\
3. An advantage of the above approach is that alternative definitions
on real axis coincide. For instance, fractional derivatives can 
be interpreted using the fractional power of the negative Laplacian
$-\Delta$, see Lemma 1 in \cite{yang10}. This can be introduced via 
Fourier transform, see Section 2.6 in \cite{samko93}. For more 
information and multidimensional extension of the Riesz derivative 
see also Section 2.10 in \cite{kilbas06}. Note that finite difference 
discretizations for Riesz fractional derivatives has been studied also
in \cite{shen08} highlighting its connection with probabilistic models.\\
4. We have left open the question for which functions does Definition 
\ref{fr_der} make sense. The general answer requires the discussion 
of the Triebel--Lizorkin spaces \cite{adams03}, \cite{samko93}, which is
beyond the scope of this paper. Some sufficient conditions on a bounded
interval $(a,b)$ are also discussed in Lemma 2.2 in \cite{kilbas06}.
At the same time, since we will approximate the Riesz derivative with 
finite differences of the fractional integrals and verify that the 
fractional integrals make sense on the function space $C_I(\er)$.

\begin{lem}
For each function $f\in C_I(\er)$ and any exponent $\alpha\in (1,2)$
the fractional order integral operators $_{-\infty}^{}I_x^{2-\alpha}f$ and 
$_{x}^{}I_\infty^{2-\alpha} f$ make sense.
\end{lem}
\emph{Proof:}
We prove the statement for the right-sided approximation, the left 
sided can be handled in a similar way. In concrete terms, we prove
that
\bee\label{right_side_finite}
\int_x^\infty \frac{f(s)}{(s-x)^{\alpha-1}} \:\mathrm{d}s < \infty.
\eee
Obviously, there is a $k\in\mathbb{Z}$ such that $a+(b-a)k > x$ and 
accordingly, 
\bee\label{right_side_finite2}
\int_x^\infty \frac{f(s)}{(s-x)^{\alpha-1}} \:\mathrm{d}s =
\int_x^{a+(b-a)k} \frac{f(s)}{(s-x)^{\alpha-1}} \:\mathrm{d}s +
\int_{a+(b-a)k}^\infty \frac{f(s)}{(s-x)^{\alpha-1}}\:\mathrm{d}s.
\eee
Here, using the condition $0<\alpha-1<1$ we have that the first term 
is finite. To estimate the second one, we introduce the
function $F: (a+(b-a)k,\infty)\to\er$ with 
$$
F(s) = \int_{a+(b-a)k}^s f(s^*) \:\mathrm{d}s^*
$$
such that $F'(s) = f(s)$ on $(a+(b-a)k,\infty)$. Also, since 
$f\in C_I(\er)$, we have that $0 = F(a+(b-a)k) = 
F(a+(b-a)(k+1)) = \dots$ such that $F$ is bounded. \\

With this the second term on the right hand side of  
\eqref{right_side_finite2} can be rewritten as
$$
\bali
\int_{a+(b-a)k}^\infty \frac{F'(s)}{(s-x)^{\alpha-1}} \:\mathrm{d}s 
&= \left[\frac{F(s)}{(s-x)^{\alpha-1}}\right]^{\infty}_{a+(b-a)k}
- \int_{a+(b-a)k}^\infty \frac{F(s)\cdot(1-\alpha)}
{(s-x)^{\alpha}} \:\mathrm{d}s \\
&= 
- \int_{a+(b-a)k}^\infty \frac{F(s)\cdot(1-\alpha)}
{(s-x)^{\alpha}} \:\mathrm{d}s,
\eali
$$ 
which is also finite since $F$ is bounded and $2-\beta>1$.
Therefore, \eqref{right_side_finite} is finite, which proves the lemma.
\quad$\square$

\paragraph{Fundamental solutions}
For the forthcoming analysis we analyze the Cauchy problem
\bee\label{basic1}
\begin{cases}
\pa_t u(t,x) = \pa_{|x|}^\alpha u(t,x) \quad t\in\erp,\; x\in\er\\
u(0,x) = u_0(x)\quad x\in\er
\end{cases}
\eee
with a given initial function $u_0\in C_I(\er)$ and $\alpha\in (1,2]$.
\begin{lem}\label{fund_sol_lem}
The Cauchy problem in \eqref{basic1} has a unique solution and 
can be given by 
\bee\label{alt_mo}
u(t,x) = (\Phi_{\alpha,t} * u_0) (x),
\eee
where $\Phi_{\alpha,t}$ denotes the fundamental solution corresponding to the 
Riesz fractional differential operator $\pa_{|x|}^\alpha$, furthermore, 
$u(t,\cdot) \in C^\infty(\er)$ for all $t\in\erp$.
\end{lem}
\emph{Proof:} 
Applying the spatial Fourier transform $\mcf$ to the 
equations in \eqref{basic1} we have  
$$
\begin{cases}
\pa_t \mcf u(t,s) = - |s|^\alpha \mcf u(t,s) \quad t\in\erp,\; s\in\er\\
\mcf u(0,s) = \mcf u_0(s)\quad s\in\er,
\end{cases}
$$
where we have used the identity 
$\mcf \left[\pa_{|x|}^\alpha u(t,x)\right] (s) = - |s|^\alpha \mcf u(t,s)$,
see 
\cite{hilfer08}, p. 38.  
Therefore, $\mcf u(t,s) = e^{-t|s|^\alpha}\mcf u_0(s)$ such that an inverse 
Fourier transform  $\mcf^{-1}$ implies that 
$$
u(t,x) = \mcf^{-1}(e^{-t|s|^\alpha} \mcf u_0(s))(x) = \mcf^{-1}(e^{-t|s|^\alpha})
* u_0 (x).
$$
In this way, the fundamental solution of \eqref{basic1} can be given as 
$$ 
\Phi_{\alpha,t} (s) = \mcf^{-1} (e^{-t|x|^\alpha})(s) = \mcf (e^{-t|x|^\alpha})(s).
$$
Using the fact that the function to transform is even and the 
Fourier transform  $\mcf \exp\{-a|\cdot|\}(s)$ can be given (see, 
\emph{e.g.}, \cite{gradstein00}, p. 1111) we have
\bee\label{f_long}
\bali
\Phi_{\alpha,t} (s) &= \mcf (e^{-t|x|^\alpha})(s) =
\mcf (e^{-t|x|} e^{-t|x|^{\alpha-1}})(s) =  
\mcf (e^{-t|x|})(s) * \mcf (e^{-t|x|^{\alpha-1}})(s) \\
&=\frac{\sqrt{\frac{2}{\pi}}t}{t^2+s^2} * \mcf (e^{-t|x|^{\alpha-1}})(s).
\eali
\eee
Here for any fixed $t$ the real function given by 
$s\to \frac{\sqrt{\frac{2}{\pi}}t}{t^2+s^2}$ is in $C^\infty(\er)$ and 
also all of its derivatives are in $L_1(\er)$. 
On the other hand, for $\alpha >1$ the real function given with 
$e^{-t|x|^{\alpha-1}}$ is in $L_1(\er)$, therefore 
$\mcf (e^{-t|x|^{\alpha-1}})$ is bounded and continuous.
Consequently, using \eqref{f_long} the right hand side of the equality 
$$
\pa^k \Phi_{\alpha,t} (s) = 
\pa^k \frac{\sqrt{\frac{2}{\pi}}t}{t^2+s^2} *  \mcf (e^{-t|x|^{\alpha-1}})(s)
$$
makes sense, which gives statement in the lemma. \quad$\square$

\paragraph{Discretization}
The finite difference approximation of the fractional order 
derivatives is not straightforward.
It turns out that an obvious one-sided finite difference 
approximation of the one-sided
Riemann--Liouville derivatives results in an unstable method even
if an implicit Euler method is applied for the time marching scheme
\cite{meerschaert04}. To stabilize these schemes, we need to use the 
\emph{translated Gr\"unwald--Letnikov formula}, which is given for
$f\in C(\er)$ with
\begin{equation}\label{trGLformula_left}
\bali
D^{\alpha,p,h}_{-\infty, GL} f(x) &= \frac{1}{\Gamma(-\alpha)} \lim_{M\to\infty}
\frac{1}{h^\alpha}
\sum_{k=0}^{M} \frac{\Gamma(k-\alpha)}{\Gamma(k+1)}  f(x-(k-p)h)\\
&=
\frac{1}{h^\alpha}
\sum_{k=0}^{\infty} g_k f(x-(k-p)h)
\eali
\end{equation}
and 
\begin{equation}\label{trGLformula_right}
\bali
D^{\alpha,p,h}_{\infty, GL} f(x) &= \frac{1}{\Gamma(-\alpha)} \lim_{M\to\infty}
\frac{1}{h^\alpha} \sum_{k=0}^{M} \frac{\Gamma(k-\alpha)}{\Gamma(k+1)}  
f(x+(k-p)h) \\
&=
\frac{1}{h^\alpha}
\sum_{k=0}^{\infty} g_k f(x+(k-p)h)
\eali
\end{equation}
depending on the translation parameter $p\in\mathbb{N}$, the order of 
the differentiation $\alpha\in (1,2]$ and the discretization parameter
$h\in\erp$, where we used the coefficients
$$
g_k = \frac{\Gamma(k-\alpha)}{\Gamma(-\alpha)\Gamma(k+1)}
= (-1)^k\binom{\alpha}{k}.
$$
The principle of the two-sided translated discretizations is depicted in 
Figure \ref{trans_fig}. 
\begin{figure}
	\setlength{\unitlength}{1.4cm}
\begin{center}
		\begin{picture}(10,1.5)(0,0)
			
			\linethickness{1pt}	
			\put(0.2,0.84){\line(1,0){9.6}}			

			\linethickness{0.75pt}	
			\put(0.5,0.62){{\Large *}}	
			\put(1.4,0.62){{\Large *}}
			\put(2.3,0.62){{\Large *}}
                        \put(3.2,0.62){{\Large *}}
			\put(4.1,0.62){{\Large *}}
			\put(5.0,0.62){{\Large *}}
                        \put(5.9,0.62){{\Large *}}
		        \put(6.8,0.62){{\Large *}}	
			 
                        \put(3.2,0.75){$\Diamond$}
			\put(4.1,0.75){$\Diamond$}	
			\put(5.0,0.75){$\Diamond$}
			\put(5.9,0.75){$\Diamond$}
                        \put(6.8,0.75){$\Diamond$}
			\put(7.7,0.75){$\Diamond$}
			\put(8.6,0.75){$\Diamond$}
		        \put(9.5,0.75){$\Diamond$}

			
		
			\put(0.4,0.4){\small{$\cdots$}}		
			\put(1.05,0.4){\small{$\displaystyle x-4h$}}
			\put(2.0,0.4){\small{$\displaystyle x-3h$}}
                	\put(2.95,0.4){\small{$\displaystyle x-2h$}}	
			\put(3.9,0.4){\small{$\displaystyle x-h$}}
			\put(5,0.4){\small{$\displaystyle x$}}
			
			\put(5.55,0.4){\small{$\displaystyle x+h$}}	
			\put(6.6,0.4){\small{$\displaystyle x+2h$}}
                        \put(7.55,0.4){\small{$\displaystyle x+3h$}}
			\put(8.5,0.4){\small{$\displaystyle x+4h$}}
		\put(9.5,0.4){\small{$\cdots$}}
		\end{picture}
		\end{center}
\caption{Basis points for the the left-sided ($*$) and the 
right-sided ($\Diamond$) translated Gr\"unwald--Letnikov
formula (given in \eqref{trGLformula_left} and \eqref{trGLformula_right})  
applied in $x$ with the translation parameter $p=2$.}
\label{trans_fig}
\end{figure}
These coefficients satisfy the following:
\bee\label{gk1}
\begin{aligned}
&\sum_{k=0}^\infty g_k = 0 \quad\forall \;\alpha\in (1,2]\\
&g_1 = -\alpha, \; g_j \ge 0\quad\textrm{for}\; j\not= 1. 
\end{aligned}
\eee
We use the same notation for the discrete differential 
(or difference) operators, i.e. for each 
$\vv=(\dots, v_{-1},v_0,v_1,\dots)\in\er^{\mathbb{Z}}$ 
we write
\bee\label{disc_gl_def}
\left[D^{\alpha,p,h}_{-\infty, GL} \vv\right]_j  =
\frac{1}{h^\alpha} \sum_{k=0}^\infty (-1)^k g_k v_{j+p-k}
\;\textrm{and}\;
\left[D^{\alpha,p,h}_{\infty, GL} \vv\right]_j   =
\frac{1}{h^\alpha} \sum_{k=0}^\infty (-1)^k g_k v_{j+k-p},
\eee
where the superscript $j$ denotes the $j$th component.\\
\emph{Remarks:}\\
1. One can prove \cite{meerschaert04} that the integrals in 
\eqref{trGLformula_left} and \eqref{trGLformula_right} approximate 
the corresponding Riemann--Liouville derivatives in the following sense:
\begin{equation}\label{trGLformula_left_approx}
_{-\infty}^{RL}\pa_x^\alpha f(x) = 
 D^{\alpha,p,h}_{-\infty, GL} f(x) + \mathcal{O}(h)
\end{equation}
and similarly, 
\begin{equation}\label{trGLformula_right_approx}
_{x}^{RL}\pa_\infty^\alpha f(x) = 
 D^{\alpha,p,h}_{\infty, GL} f(x) + \mathcal{O}(h),
\end{equation}
provided that both of the Fourier transform of $f$ and that of
$_{-\infty}^{RL}\pa_x^\alpha f(x)$ and $_{x}^{RL}\pa_\infty^\alpha f(x)$ 
are in $L_1(\er)$. We will point out that in the present framework 
no such assumption is necessary.\\
2. Higher-order finite difference approximations can be obtained 
as a linear combination of first order ones using different translation 
parameters. 
For example, the sum defined by 
\bee\label{2ndorder_RL}
_{-\infty}^{GL}D_{x,h}^{\alpha,p,q}u(x)=
\frac{2q-\alpha}{2(p-q)}
D^{\alpha,p,h}_{-\infty, GL} u(x)
+
\frac{2p-\alpha}{2(p-q)}
D^{\alpha,p,h}_{-\infty, GL} u(x)
\eee
provides a second order accurate approximation of the Riemann--Liouville
derivative.
Similar statement holds if $-\infty$ is switched to $\infty$.
For the details, see \cite{zhou13}.\\ 
3. The nonlocal effect of the differential operators result in full
matrices. At the same time, one can save some computing efforts with 
an appropriate decomposition of the above matrix \cite{wang12}.\\
4. An alternative approximation of fractional order elliptic operators
(which can be applied in multidimensional cases) was proposed in 
\cite{bonito13}, which can be a basis also for finite element 
discretizations.






\section{Results}

\subsection {Extensions }
Extensions are not only necessary to have well-posed problems involving 
nonlocal diffusion operators, but also essential at the discrete level.  
In order to have sufficient accuracy in the finite difference 
approximation near to the boundary and at the boundary of a nonlocal
differential operator, it is necessary to have (virtual) gridpoints 
outside of the original computational domain $\Omega = (a,b)$.  
This is clearly shown in \eqref{trGLformula_left}, 
\eqref{trGLformula_right} and \eqref{2ndorder_RL}.

Summarized, the sketch of our approach is the following:
\begin{itemize}
\item 
we extend the problem to $\er$ to get rid of the boundary conditions
\item
we solve the corresponding Cauchy problem  \eqref{basic1}
(we will approximate this with finite differences)
\item
we verify that the desired homogeneous Neumann (or homogeneous Dirichlet)
boundary conditions are satisfied for the restriction of the solution. 
\end{itemize}


\begin{dfn}\label{compatible_def}
We say that the extension 
\[
\tilde\cdot: L_2(\Omega) \to C_I(\er)
\]
is \emph{compatible} with the homogeneous Neumann (no-flux) or
Dirichlet boundary conditions and the operator $\pa_{|x|}^\alpha$ if the 
function $\tilde u$ is the unique solution of the following problem
\[
\begin{cases}
\pa_t \tilde u (t, x) = \pa_{|x|}^\alpha \tilde u(t,x)\quad 
t\in\erp, \; x\in\er\\
\tilde u (0, x) = \tilde u_0(x)\quad x\in\er\\
\end{cases}
\]
and $\pa_x \tilde u (t, x) = 0$  or $\tilde u (t, x) = 0$ for 
$x\in\pa\Omega$ and $t>0$.
\end{dfn}
In rough terms, one can say that a correct extension of the solution
from $\Omega$ is the function which solves the extended problem on 
the whole real axis such that the boundary condition on $\pa\Omega$ 
is still satisfied. 


\paragraph{Extension corresponding to homogeneous Dirichlet
boundary condition}
As a motivation, we use the idea in \cite{szymczak03} which is generalized 
to the case of two absorbing walls. For the simplicity, the definition is 
given for functions $u:(0,1)\to\er$.
\begin{dfn}
We call the 2-periodic extension of the function
\bees
\hat u^{\mcd,1} (x) = 
\begin{cases}
u(x)\quad\textrm{for}\; x\in (0,1)\\
-u(-x) \quad\textrm{for}\; x\in (-1,0)
\end{cases}
\eees
the Dirichlet type extension of $u$ and we denote it 
with $\hat u^{\mcd}$.
\end{dfn}
\noindent\emph{Remarks:}\\
1. This is sometimes called the odd extension of $u$ and can be 
obtained first with reflecting the graph of $u$ to 
$(0, 0)$ and $(1, 0) \in\er^2$ to extend it to $(-1, 2)$ and 
reflecting this graph further  to $(-1, 0)$ and $(2, 0) \in\er^2$ to 
extend it to $(-2, 3)$ and continuing this process. \\
2. A natural physical interpretation of this extension is the following:
To ensure zero-concentration at the end points in the consecutive 
time steps, we have to force a skew-symmetric concentration profile 
around $0$ and $1$.\\ 
3. We may define $u^\mcd (k) =0$ for $k\in\mathbb{Z}$ but as we will
see this is not essential.

A simple calculation shows that we have
\bee\label{eq_for_hatD2}
\hat u^{\mcd} (y) = -\hat u^{\mcd} (-y)\quad\textrm{and}\quad
\hat u^{\mcd} (1+y) = -\hat u^{\mcd} (1-y)\quad y\in\er. 
\eee

We define similarly the extension $\hat \vv^{\mcd}\in \er^{\mathbb{Z}}$ of 
the vector $\vv = (v_0, v_1, \dots, v_{N+1}) \in\er^{N+2}$ with 
$v_0 = v_{N+1} =0$:
\bee\label{ext_vector}
[\hat \vv^{\mcd}]_j = 
\begin{cases}
v_j \quad j = 0, 1, \dots, N+1\\ 
-v_{2(N+1)-j} \quad j = N+2, N+3, \dots, 2N+1\\
v_{j+2N} \quad j \not\in\left\{ 0, 1, \dots, 2N+1\right\}.
\end{cases}
\eee


\paragraph{Extension corresponding to homogeneous Neumann
boundary condition}
Similarly to the previous case, the definition is 
given for functions $u:(0,1)\to\er$.

\begin{dfn}
We call the 2-periodic extension of the function
\bees
\hat u^{\mcn,1} (x) = 
\begin{cases}
u(x) \quad\textrm{for}\; x\in (0,1)\\
u(-x) \quad\textrm{for}\; x\in (-1,0)
\end{cases}
\eees
the Neumann type extension of $u$ and we denote it 
with $\hat u^{\mcn}$.
\end{dfn}
\noindent\emph{Remarks:}\\
1. This is sometimes called the even extension of $u$ and can be 
obtained first with reflecting the graph of $u$ to the vertical 
line given with $x=0$ to extend it to $(-1,0)$ then to the vertical 
line given with $x=1$ to extend it to $(1,2)$ and continuing this 
process.\\
2. A natural physical interpretation of this extension is the following:
To ensure zero-flux, we force a symmetric concentration profile 
around $0$ and $1$. 

A simple calculation shows that 
\bee\label{eq_for_hatN2}
\pa_x \hat f^{\mcn} (1+y) = \pa_x \hat f^{\mcn} (1-y)\quad 
\textrm{and} \quad
\pa_x \hat f^{\mcn} (y) = \pa_x \hat f^{\mcn} (-y)\quad y\in\er. 
\eee
Similar notations are used for the ``extended'' vector 
$\hat \vv^\mcn\in \er^{\mathbb{Z}}$ of
$\vv=(v_0, v_1, \dots, v_{N})\in\er^{N+1}$ which is defined as follows:
\bee\label{ext_vector_n}
[\hat \vv^{\mcn}]_j =
\begin{cases}
v_j \quad j = 0, 1, \dots, N\\
v_{-j-1} \quad j = -N-1, -N, \dots, -1\\
v_{j+2N+2} \quad j \not\in\left\{ -N-1, \dots, 0,\dots, N\right\}.
\end{cases}
\eee
A natural physical interpretation of this extension is that
particles are reflected at the boundaries to ensure zero flux.
In this way, we also reflect the concentration profile in the model.
The principle of the Neumann extension for a vector is visualized in
Figure \ref{neumann_ext_fig}.
\begin{figure}
	\setlength{\unitlength}{1.4cm}
\begin{center}
		\begin{picture}(10,1.5)(0,0)
			
			\linethickness{1pt}	
			\put(0.2,0.85){\line(1,0){9.6}}			

			\linethickness{0.75pt}	
			\put(0.4,1.10){$\diamond$}	
			\put(1.2,1.60){$\diamond$}
			\put(2.0,1.35){$\diamond$}
                        \put(2.8,1.85){$\diamond$}
			\put(3.6,1.85){$\bullet$}
			\put(4.4,1.35){$\bullet$}
                        \put(5.2,1.60){$\bullet$}
		        \put(6.0,1.10){$\bullet$}	
                        \put(6.8,1.10){$\diamond$}
			\put(7.6,1.60){$\diamond$}
			\put(8.4,1.35){$\diamond$}
                        \put(9.2,1.85){$\diamond$}

                        \put(0.4,0.6){$v_{-4}$}	
			\put(1.2,0.6){$v_{-3}$}
			\put(2.0,0.6){$v_{-2}$}	
                        \put(2.8,0.6){$v_{-1}$}	
			\put(3.6,0.6){$v_{0}$}	
			\put(4.4,0.6){$v_{1}$}	
                        \put(5.2,0.6){$v_{2}$}	
		        \put(6.0,0.6){$v_{3}$}		
                        \put(6.8,0.6){$v_{4}$}	
			\put(7.6,0.6){$v_{5}$}	
			\put(8.4,0.6){$v_{6}$}	
                        \put(9.2,0.6){$v_{7}$}	

		\end{picture}
		\end{center}
\caption{Entries of the Neumann extension $\hat \vv^\mcn$ of the vector 
$\vv = (v_0, v_1, v_2, v_3)$.} 
\label{neumann_ext_fig}
\end{figure}
We verify that the above extensions meet the requirements in 
Definition \ref{compatible_def}.

\begin{lem}\label{compatible_lem}
The extensions $\hat u^\mcn$ and $\hat u^\mcd$ 
are compatible with the Dirichlet and the Neumann boundary
conditions, respectively and with the differential operator
$\pa_{|x|}^\alpha$. 
\end{lem}
\emph{Proof:} 
According to Lemma \eqref{fund_sol_lem}
we can express the solution of 
\[
\begin{cases}
\pa_t u = \pa_{|x|}^\alpha u\quad \textrm{on}\; \erp\times\er\\
u(0,\cdot) = \hat u^\mcn
\end{cases}
\]
with the convolution 
\[
u(t, \cdot) = \Phi_{\alpha,t} * \hat u^\mcn = \hat u^\mcn * \Phi_{\alpha,t},   
\]
such that 
\[
u(t,x) = \int_\er \hat u^\mcn(x-y) \Phi_{\alpha, t} (y)\dy.  
\]
Since $\Phi_{\alpha, t}\in C^\infty(\er)$, the same holds for $u(t,\cdot)$.
Accordingly, the right and left limit of $\pa_x u(t,\cdot)$ in $1$ coincide.
Using this, \eqref{eq_for_hatN2} and the fact that $\Phi_{\alpha, t}$ is even
we obtain  
\[
\bali
\pa_x u(t,1) &= \lim_{\epsilon_{n}\to 0-} \pa_x u(t, 1-\epsilon_n)
= 
\lim_{\epsilon_n\to 0-} (\pa_x \hat u^\mcn * \Phi_{\alpha,t})(1-\epsilon_n)\\
&=
\lim_{\epsilon_n\to 0-} \int_\er \pa_x\hat u^\mcn(1-\epsilon_n - y) 
\Phi_{\alpha, t} (y)\dy 
=
-\lim_{\epsilon_n\to 0-} \int_\er \pa_x\hat u^\mcn(1+\epsilon_n + y) 
\Phi_{\alpha, t} (y)\dy\\ 
&=
-\lim_{\epsilon_n\to 0-} \int_\er \pa_x\hat u^\mcn(1+\epsilon_n + y) 
\Phi_{\alpha, t} (-y)\dy\\
&= 
-\lim_{\epsilon_n\to 0-} \int_\er \pa_x\hat u^\mcn(1+\epsilon_n - y) 
\Phi_{\alpha, t} (y)\dy
= 
-\lim_{\epsilon_n\to 0-} (\pa_x \hat u^\mcn * \Phi_{\alpha,t})(1+\epsilon_n)\\
&=
-\lim_{\epsilon_n\to 0-} \pa_x u(t, 1+\epsilon_n) = - \pa_x u(t,1),
\eali
\]
which gives that the homogeneous Neumann boundary condition is satisfied
in $1$. With an obvious modification, using \eqref{eq_for_hatN2}, we can 
also verify the homogeneous Neumann boundary condition $\pa_x u(t,0)=0$.

Similarly, we can express the solution of 
\[
\begin{cases}
\pa_t u = \pa_{|x|}^\alpha u\quad \textrm{on}\; \erp\times\er\\
u(0,\cdot) = \hat u^\mcd
\end{cases}
\]
with 
\[
u(t,x) = \int_\er \hat u^\mcd(x-y) \Phi_{\alpha, t} (y)\dy.  
\]
Since $\Phi_{\alpha, t}\in C^\infty(\er)$, the same holds for $u(t,\cdot)$.
Accordingly, the right and left limit of $u(t,\cdot)$ in $1$ coincide.
Using this and \eqref{eq_for_hatD2} we obtain  
\[
\bali
u(t,1) &= \lim_{\epsilon_{n}\to 0-}  u(t, 1-\epsilon_n)
= 
\lim_{\epsilon_n\to 0-}  \hat u^\mcd * \Phi_{\alpha,t} (1-\epsilon_n)\\
&=
\lim_{\epsilon_n\to 0-} \int_\er \hat u^\mcd (1-\epsilon_n - y) 
\Phi_{\alpha, t} (y)\dy 
=
-\lim_{\epsilon_n\to 0-} \int_\er \hat u^\mcd (1+\epsilon_n + y) 
\Phi_{\alpha, t} (y)\dy\\ 
&=
-\lim_{\epsilon_n\to 0-} \int_\er \hat u^\mcd (1+\epsilon_n + y) 
\Phi_{\alpha, t} (-y)\dy
= 
-\lim_{\epsilon_n\to 0-} \int_\er \hat u^\mcd (1+\epsilon_n - y) 
\Phi_{\alpha, t} (y)\dy\\ 
&= 
-\lim_{\epsilon_n\to 0-} \hat u^\mcd  * \Phi_{\alpha,t}(1+\epsilon_n)
=
-\lim_{\epsilon_n\to 0-}  u(t, 1+\epsilon_n) = - u(t,1),
\eali
\]
which gives that the homogeneous Dirichlet boundary condition is satisfied
in $1$. With an obvious modification, using \eqref{eq_for_hatD2}, we can 
verify homogeneous Dirichlet boundary condition also in $0$.

In both cases, the equality $\pa_t u = \pa_{|x|}^\alpha u$ is satisfied 
in $(0, 1)$ which gives the statement in the lemma.\quad $\square$

Using the above extension, we will investigate the numerical solution of 
the problem
\bee\label{ext_PDE_N2}
\begin{cases}
\pa_t u(t,x)  = 
\pa_{|x|}^\alpha u^\mcn (t,\cdot)(x)
\quad t\in (0,T),\; x\in (0,1)\\
u(0,x) = u_0 \quad  x\in (0,1).  
\end{cases}
\eee

The following theorem is the basis of our numerical method, which again 
confirms the favor of the Neumann type extension. 
\begin{thm}\label{compatible_thm}
For any $u_0\in C[0,1]$
the problem in \eqref{ext_PDE_N2} is well-posed, and its unique 
solution $u\in C^\infty [0,1]$
satisfies the boundary conditions $\pa_x u(t,0) = \pa_x u(t,1) =0$.
\end{thm}
\emph{Proof:}
We first note that the problem
\bee\label{ext_PDE_N2_full}
\begin{cases}
\pa_t u(t,x)  = \pa_{|x|}^\alpha  u(t,x) \quad t\in (0,T),\; x\in\er\\
u(0,x)  = \hat u_0^\mcn (x)\quad  x\in\er.  
\end{cases}
\eee
is well-posed and corresponding to Lemma \ref{compatible_lem} its 
solution can be given by
\[
u(t,x) = \int_\er \hat u^\mcn(x-y) \Phi_{\alpha, t} (y)\dy.  
\]
This implies that
\[
\bali
u(t,-x) 
&= \int_\er \hat u^\mcn(-x-y) \Phi_{\alpha, t} (y)\dy =
\int_\er \hat u^\mcn(x+y) \Phi_{\alpha, t} (-y)\dy\\
&=
\int_\er \hat u^\mcn(x-y) \Phi_{\alpha, t} (y)\dy = u(t,x).
\eali  
\]
and 
\[
\bali
u(t,x+2) 
&= \int_\er \hat u^\mcn(x+2-y) \Phi_{\alpha, t} (y)\dy =
\int_\er \hat u^\mcn(x-y) \Phi_{\alpha, t} (y)\dy = u(t,x).
\eali  
\]
Therefore, the restriction of the solution $u$ of \eqref{ext_PDE_N2_full}
solves the problem in \eqref{ext_PDE_N2}.
Here we have also used throughout that for the Neumann extension: 
$u^\mcn\in C_I(\er)$ 
such that the Riesz derivative $\pa^{\alpha}_{|x|}u^\mcn$ makes sense.

To prove the uniqueness, we assume that $v$ solves \eqref{ext_PDE_N2}. 
According to \eqref{eq_for_hatN2} we obviously have that the Neumann 
extension satisfies
\bee\label{t_der_kit}
\pa_t \hat v^\mcn (t,-x) = \pa_t \hat v^\mcn (t,x) \quad x\in\er.
\eee
To compute the fractional order integral we introduce the function 
$C_I(\er)\ni \hat v_0^\mcn = \hat v(t,\cdot)^\mcn - \int_0^1 v(t,\cdot)$ 
such that we have 
$$
\bali
&\Gamma(2-\alpha)
_{\infty}I_x^{2-\alpha} \hat v^\mcn_0 (-x)
=
\int_{-\infty}^{-x} \frac{v^\mcn_0(y^*)}{(-x-y^*)^\alpha}\mathrm{d}y^* 
+ 
\int^{\infty}_{-x} \frac{v^\mcn_0(y^*)}{(y^*+x)^\alpha}\mathrm{d}y^*\\
&=  
\int^{\infty}_{x} \frac{v^\mcn_0(-y)}{(-x+y)^\alpha}\dy 
+ 
\int_{-\infty}^{x} \frac{v^\mcn_0(-y)}{(-y+x)^\alpha}\dy
=  
\int^{\infty}_{x} \frac{v^\mcn_0(y)}{(-x+y)^\alpha}\dy 
+ 
\int_{-\infty}^{x} \frac{v^\mcn_0(y)}{(-y+x)^\alpha}\dy\\
&=
\Gamma(2-\alpha)
_{\infty}I_x^{2-\alpha} \hat v^\mcn_0 (x).
\eali
$$
Therefore, we also have 
\bee\label{xder_neumannext}
\pa^\alpha_{|x|} v^\mcn(t,-x) = \pa^\alpha_{|x|} v^\mcn(t,x).
\eee
Consequently, \eqref{t_der_kit} and \eqref{xder_neumannext} imply
$$
 \pa_t \hat v^\mcn (t,-x) = \pa_t \hat v^\mcn (t,x) =
\pa^\alpha_{|x|} v^\mcn(t,x) = \pa^\alpha_{|x|} v^\mcn(t,-x)
$$
and the periodicity obviously gives 
$$
 \pa_t \hat v^\mcn (t,2+x) = \pa_t \hat v^\mcn (t,x) =
\pa^\alpha_{|x|} v^\mcn(t,x) = \pa^\alpha_{|x|} v^\mcn(t,2+x)
$$
such that the Neumann extension $u^\mcn$ solves \eqref{ext_PDE_N2_full}.
This, however, has a unique solution,
which gives the uniqueness of the solution of \eqref{ext_PDE_N2}.
\quad$\square$

\subsection{Numerical methods}
Following the classical method of lines technique we first 
discretize the spatial variables in the extended problem 
and choose a time stepping scheme for the full discretization.




To streamline the forthcoming computations, the interval $[0, 1]$ will 
be transformed to $[0, \pi]$, where we use the following grid points:
\bee\label{ezxj}
x_j:= \pi h\left(j + \frac{1}{2}\right) = 
\pi\left(\frac{1}{2(N+1)} + \frac{j}{N+1}\right).
\eee
 
$u^n_j$ denotes the numerical approximation at time $n\tau$ in the grid point 
$x_j$ and $\uu(n\delta,\cdot)$ the values of the analytic solution at time 
$n\delta$ in the grid points.

\subsubsection{Analysis of a finite difference scheme}

For the spatial discretization we use the Gr\"unwald--Letnikov 
approximations in \eqref{disc_gl_def} introducing 
$A_{\alpha,h}\in \er^{(N+1)\times(N+1)}$ with 
\bee\label{A_def}
A_{\alpha,h} \uu = D_{-\infty,\textrm{GL}}^{\alpha,1,h} \uu^{\mcn}
         + D_{\infty,\textrm{GL}}^{\alpha,1,h} \uu^{\mcn}.
\eee
This is combined with an implicit Euler time stepping to obtain
\bee\label{impl_euler}
\frac{u_j^{n+1} - u_j^{n}} {\tau} 
= 
\left[A_{\alpha,h} u^{n+1}\right]_j.
\eee
To make the consecutive formulas more accessible, we 
expand \eqref{A_def} in a concrete example.\\
\textbf{Example} \emph{We give the first component of 
$A_{\alpha,h}\vv$ for $\vv=(v_0, v_1, v_2, v_3)$.}
$$
\bali
\left[A_{\alpha,h}\vv\right]_1
& = g_0 v_0 + g_1 v_1 + g_2 v_2 + g_3 v_3 + g_4 v_3 + g_5 v_2 
+ g_6 v_1 + g_7 v_0 + g_8 v_0 + \dots \\
& + g_0 v_2 + g_1 v_1 + g_2 v_0 + g_3 v_0 + g_4 v_1 + g_5 v_2 
+ g_6 v_3 + g_7 v_3 + g_8 v_2 + \dots.
\eali
$$
In the general case, using \eqref{A_def}, \eqref{disc_gl_def}
and \eqref{ext_vector_n} we obtain that 
\bee\label{give_A1}
\bali
& [A_{\alpha,h}\vv]_j = \frac{C_\sigma}{h^\alpha}
\left(
\sum_{k=0}^{j+1} g_k v_{j-k+1} + \sum_{k=0}^{N-j+1} g_k v_{k+j-1} 
\right.\\
&
\left. +
\sum_{l=0}^{\infty} \left(\sum_{k=j+2}^{j+N+2} g_{k+2l(N+1)} v_{k-j-2}
+\sum_{k=N-j+2}^{2N-j+2} g_{k+2l(N+1)} v_{2N-k-j+2}\right)\right.\\
&\left. +
\sum_{l=0}^{\infty} \left(\sum_{k=j+N+3}^{j+2N+3}g_{k+2l(N+1)} v_{2N+j-k+3}
+\sum_{k=2N-j+3}^{3N-j+3}g_{k+2l(N+1)} v_{k-2N+j-3}
\right)\right)\\
& j=1,2,\dots,N-1
\eali
\eee
\bee\label{give_A2}
\bali
& [A\vv]_0= \frac{2C_\sigma}{h^\alpha}
\left(g_0 v_1 + g_1 v_0\right.\\
& \left. +
\sum_{l=0}^{\infty} \left(\sum_{k=2}^{N+2} g_{k+2l(N+1)} v_{N+2-k}
                         +\sum_{k=N+3}^{2N+3} g_{k+2l(N+1)} v_{2N-k+3}
\right)\right),
\eali
\eee
\bee\label{give_A3}
\bali
&[A\vv]_N = \frac{2C_\sigma}{h^\alpha}
\left(
g_0 v_{N-1} + g_1 v_N\right.\\ 
& \left. +
\sum_{l=0}^{\infty} \left(\sum_{k=2}^{N+2} g_{k+2l(N+1)} v_{N-k+2}
                         +\sum_{k=N+3}^{2N+3} g_{k+2l(N+1)} v_{k-N-3}
\right)\right).
\eali
\eee
For $K=2m+1$ the corresponding matrix can be given with a slight 
modification. Observe that all of the coefficients $g_i$ arises once on 
the right hand side of \eqref{give_A1}, \eqref{give_A2} and \eqref{give_A3}.

\begin{prop}\label{prop_matrix}
The matrix $A_{\alpha,h}$ has negative diagonal and non-negative 
off-diagonal elements.
\end{prop}
\emph{Proof:}
Observe that in \eqref{give_A1}, \eqref{give_A2} and \eqref{give_A3} 
the coefficient of $v_j, v_0$ and $v_N$, respectively is $g_1$ and $g_1$
appears only here. Therefore, using also \eqref{gk1} we obtain that 
$A$ has negative diagonals and positive off-diagonals.\quad$\square$


We analyze the properties of the matrix  $A_{\alpha,h}$ and the corresponding 
differential operator. 
\begin{lem}\label{lem_eig}
The eigenvectors of the matrix 
$A_{\alpha,h}\in \er^{(N+1)\times(N+1)}$ are given as 
\[
\vv_h^k = (\cos kx_0 , \cos kx_1, \dots, \cos kx_{N-1}, \cos kx_N)^T
\]
for each $k=0, 1, 2 \dots, N$ with the corresponding eigenvalues
\[
-\frac{\sigma}{\cos \left(\frac{\alpha}{2}\pi\right)}
\left(\frac{2}{h}\right)^\alpha
\sin^\alpha \frac{k\pi h}{2} \cdot 
\cos \left(k\pi h + \frac{\alpha}{2}(\pi-k\pi h)\right).
\]
\end{lem}

\begin{lem}\label{eig_cos_lem}
The eigenfunctions of the operator $\pa^{\alpha}_{|x|}$ on $(0,1)$ 
with homogeneous Neumann boundary conditions are given by
$\left\{ \cos k\pi x\right\}_{k=1}^{\infty}$ with the corresponding 
eigenvalues $\left\{\sigma\cdot(k\pi)^\alpha\right\}_{k=1}^{\infty}$.
\end{lem}
The technical proofs of these statements are postponed to the Appendix.

\begin{prop}\label{prop_cons}
The numerical approximation defined in the scheme \eqref{impl_euler} 
is consistent with the problem in \eqref{ext_PDE_N2} in the maximum norm
and the order of the consistency is $\mathcal{O}(\tau) + \mathcal{O}(h)$.
\end{prop}
\emph{Proof:}
It is sufficient to prove that the right hand side of \eqref{impl_euler} 
provides
a first order  approximation for the Riesz derivative of order $|\alpha|$. 
According to the proof of Theorem \ref{compatible_thm} the analytic solution 
$u^\mcn$ of \eqref{ext_PDE_N2_full} is periodic, smooth and it satisfies the 
homogeneous Neumann boundary conditions in $0$ and $1$. Therefore, its cosine 
Fourier series is pointwise convergent:
\bee\label{fexpan}
u(t,x) = \sum_{k=0}^{\infty} F_k \cos k\pi x\quad x\in(0,1),\; t\in\erp,
\eee
where we do not denote the time dependence of the cosine Fourier 
coefficients $F_k$. 
Using again that $u(t,\cdot) \in C^\infty(\er)$ we also have - using the 
regularity theory of Fourier series - that
for any $r\in\mathbb{N}$ there exists a constant $C_r$ such that 
\bee\label{est_fcoeff}
|F_k| \le \frac {C_r}{(k\pi)^r} \quad\forall k\in\mathbb{N}^+.
\eee
Using \eqref{fexpan} componentwise for $t=n\delta$ we have that
\bee\label{comp_der_vect}
\bali
\uu(h\delta,\cdot) &= 
[u(n\delta, \frac{x_0}{\pi})\; u(n\delta, \frac{x_1}{\pi})\; \dots\;  
u(n\delta, \frac{x_N}{\pi})]^T\\
&= 
[\sum_{k=0}^\infty F_k \cos k x_0\;  \sum_{k=0}^\infty F_k \cos k x_1\;
\dots\;  \sum_{k=0}^\infty F_k \cos k x_N]^T
=  
\sum_{k=0}^\infty F_k \vv_h^k.
\eali
\eee
Using Lemma \ref{eig_cos_lem} for the expansion in \eqref{fexpan} and 
the matrix $A_{\alpha,h}$ for \eqref{comp_der_vect} according to 
\eqref{long_eig} we obtain the following equality:
\bee\label{konz1}
\bali
&
\pa_{|x|}^{\alpha} u (n\delta, \frac{x_j}{\pi}) - 
\left[A_{\alpha,h} \uu(n\delta,\cdot)\right]_j \\
&=
\sum_{k=0}^{\infty} -\sigma (k\pi)^\alpha F_k \cos k\pi x_j -
\left[A_{\alpha,h} \sum_{k=0}^\infty F_k\vv_h^k \right]_j\\
&=
\sum_{k=0}^{\infty} F_k \cos k\pi x_j 
\left(-\sigma (k\pi)^\alpha
- \frac{2C_\sigma}{h^\alpha}
2^\alpha\sin^\alpha \frac{k\pi h}{2}\cdot 
\cos \left(k\pi h + \frac{\alpha}{2}(\pi-k\pi h)\right)\right)\\
&=
-\sum_{k=0}^{\infty}\sigma F_k \cos k\pi x_j (k\pi)^\alpha
\left(1 -
\left(\frac{\sin \frac{k\pi h}{2}}{\frac{k\pi h}{2}}\right)^\alpha 
\frac{\cos \left(k\pi h + \frac{\alpha}{2}(\pi-k\pi h)\right)}
{\cos\frac{k\pi h}{2}}\right).
\eali
\eee
To prove the proposition we first verify that
for a mesh-independent constant $C$ the following inequality
is valid:
\bee\label{verify_prop}
\left|
1-\left(\frac{\sin s}{s}\right)^\alpha
\left(\cos s(\alpha-2) - \sin s(\alpha -2) \tan \frac{\alpha\pi}{2}\right)
\right|
\le
C s,
\eee
which will be applied with $s = \frac{k\pi h}{2}$.
We first verify that
\bee\label{verify_first}
h(s):=
1-\left(\frac{\sin s}{s}\right)^\alpha
\left(\cos s(\alpha-2) - \sin s(\alpha -2) \tan \frac{\alpha\pi}{2}\right)
\le
C s,
\eee
where $\lim_{s\to 0+}h(s)=0$. Therefore, it is sufficient to prove that
$h'$ is bounded on $[0, \frac{\pi}{2}]$. Obviously,
$$
\bali
h'(s) &= \alpha\left(\frac{\sin s}{s}\right)^{\alpha-1} 
\frac{\sin s - s\cdot\cos s}{s^2}\cdot (\alpha - 2)
\left[\cos s(\alpha-2) - \sin s(\alpha-2)\tan \frac{\alpha\pi}{2}\right]\\
&+
\left(\frac{\sin s}{s}\right)^\alpha
\left((\alpha-2)\left(\sin s(\alpha-2) - \cos s(\alpha -2) \tan \frac{\alpha\pi}{2}
\right)\right), 
\eali
$$
where 
$$
\lim_{0+} \frac{s\cdot\cos s - \sin s}{s^2} = 
\lim_{0+} \frac{s \sin s}{ 2s} = 0. 
$$
Hence, all components in the expansion of $h'(s)$ are bounded, which
really verifies \eqref{verify_first}.
Using \eqref{verify_first} in \eqref{konz1} and applying 
\eqref{est_fcoeff} we obtain the following estimation:
\bee\label{konz1_continued}
\bali
&
\left|\pa_{|x|}^{\alpha} u (n\delta, \frac{x_j}{\pi}) - 
\left[A_{\alpha,h} \uu(n\delta,\cdot)\right]_j\right|
\le
\frac{C}{2} h 
\left| 
\sum_{k=0}^{\infty}\sigma F_k \cos k\pi x_j (k\pi)^\alpha k\pi
\right|\\
&
\le
\frac{C\cdot\sigma}{2} h 
\sum_{k=0}^{\infty} |F_k| (k\pi)^{\alpha+1}
\le
\frac{C_5 C\cdot\sigma}{2} h 
\sum_{k=1}^{\infty} \frac{1}{(k\pi)^{5}} (k\pi)^{\alpha+1}
\le
\frac{C_5 C\cdot\sigma}{2} h 
\sum_{k=1}^{\infty} \frac{1}{(k\pi)^{2}},
\eali
\eee
which completes the proof. \quad$\square$
\begin{thm}\label{main_thm}
The numerical approximation defined in the scheme \eqref{impl_euler} 
converges to the solution of \eqref{ext_PDE_N2} in the maximum norm
for $\alpha\in (1, 2]$ and the order of the convergence 
is $\mathcal{O}(\tau) + \mathcal{O}(h)$.
\end{thm}
\emph{Proof:}
Using Proposition \ref{prop_cons} we only have to verify the stability of 
\eqref{impl_euler}. 
For this, we rewrite it into a linear system
\[
\left(I - \tau A_{\alpha,h}\right) \uu^{n+1} = \uu^n.  
\]
Using Proposition \ref{prop_matrix} we obtain that the diagonal 
of $I - \tau A_{\alpha,h}$ 
is strictly positive and has non-positive off-diagonals. Moreover, 
using \eqref{give_A1}, a simple calculation shows that for the indices 
$j= 1,2,\dots, N-1$ we have 
\bee\label{poz_matr}
\bali
&\left[\left(I - \tau A_{\alpha,h}\right)\cdot 
(1,1,\dots,1)^T\right]_j \\
&=
1 + \frac{\tau}{h} \left(\alpha - \sum_{\substack{k=0\\ k\not= 1}}^{j+1} g_k -
\sum_{l=0}^{\infty} \sum_{k=j+2}^{j+2+N} g_{k+2l(N+1)} -
\sum_{l=0}^{\infty} \sum_{k=j+N+3}^{j+2N+3}g_{k+2l(N+1)}\right)\\
&+
 \frac{\tau}{h} \left(\alpha - \sum_{\substack{k=0\\ k\not= 1}}^{N-j+1} g_k - 
\sum_{l=0}^{\infty} \sum_{k=N-j+2}^{2N-j+2} g_{k+2l(N+1)} -
\sum_{l=0}^{\infty} \sum_{k=2N-j+3}^{3N-j+3}g_{k+2l(N+1)}\right).
\eali
\eee
Observe that in the brackets in \eqref{poz_matr} each coefficient
$g_i,\; i\not =1$ appears once (see the Example and the remark after 
\eqref{give_A3}). 
According to \eqref{gk1}, we obtain
\[
\alpha - \sum_{\substack{k=0\\ k\not= 1}}^{\infty} g_k = 
- \sum_{k=0}^{\infty} g_k = 0
\]
and therefore, the sum in both brackets 
on the right hand side of \eqref{poz_matr} is zero
such that the entire right hand side is \emph{one}. 
Using \eqref{give_A2} and \eqref{give_A3}, an obvious 
modification of \eqref{poz_matr} gives its positivity both for the
indices $j=0$ and $j=N$. 
Therefore,
$$
\left[\left(I - \tau A_{\alpha,h}\right)\cdot 
(1,1,\dots,1)^T\right]_j = 1
$$
is valid for all $j=0,1,\dots, N$.
In this way, $(I-\tau A_{\alpha,h})^{-1}$ elementwise positive 
such that $\|(I-\tau A_{\alpha,h})^{-1}\|_\infty= 1$, and 
consequently, the scheme in \eqref{impl_euler} is unconditionally stable. 
\quad$\square$

\paragraph{Construction of the matrix  $A_{\alpha,h}$}
Whenever the coefficient in the matrix  $A_{\alpha,h}$
are based on an infinite number of grid points, it can be 
computed in concrete terms. 

For this we introduce $B_{\alpha,h}\in\er^{(N+1)\times (N+1)}$ with
$$
B_{\alpha,h} = \left (\vv^0_h, \vv^1_h, \dots,\vv^N_h\right), 
$$
which consists of the eigenvectors of $A_{\alpha,h}$, see 
Lemma \ref{lem_eig}. Then 
$$
(I-\tau A_{\alpha,h}) B_{\alpha,h} = 
\left((1-\tau\lambda_1)\vv^0_h, 
(1-\tau\lambda_1)\vv^1_h, \dots, (1-\tau\lambda_1)\vv^N_h\right), 
$$
and therefore,
\bee\label{give_Aalpha}
\tau A_{\alpha,h} = 
I - \left((1-\tau\lambda_1)\vv^0_h, 
(1-\tau\lambda_1)\vv^1_h, \dots, (1-\tau\lambda_1)\vv^N_h\right)
B_{\alpha,h}^{-1}, 
\eee
where on the right hand side all terms can be computed.

\paragraph{Complexity and extension to higher-order methods}
Since we have explored the eigenvectors of $A_{\alpha,h}$ we 
do not have to compute its components in the practice as a series.   
We simply obtain the matrix using \eqref{give_Aalpha} such that
the extension does not result in extra computational costs. 

Higher order methods can be obtained in the same fashion. In the 
semidiscretization, we should then choose a higher order spatial
approximation, \emph{e.g.}, the one in \eqref{2ndorder_RL} and 
the accuracy of the time stepping can also be increased, \emph{e.g.}, 
using a Crank--Nicolson scheme. For homogeneous Dirichlet boundary 
conditions, such a study is performed (even for the multidimensional 
case) in \cite{deng14}.

\section{Numerical experiments}\label{section_numex}

\paragraph{A homogeneous model problem}
We first investigate the model problem
\bees
\begin{cases}
\pa_t u(t,x) = 
0.25 \pa_{|x|}^{1.2} u(t,x)
\quad x\in (0,1),\; t\in (0, 1) \\
u(0,x) = \frac{x^4}{4} - \frac{x^2}{2}\quad x\in (0,1)\\
\pa_x u (t,0) = \pa_x u (t,1) = 0\quad t\in (0, 1),
\end{cases}
\eees
which is converted to the well-posed extended problem
\bee\label{model_pr_ext}
\begin{cases}
\pa_t u(t,x) = 
0.25 \pa_{|x|}^{1.2} u(t,x)
\quad x\in (0,1),\; t\in (0, 1) \\
u(0,x) = \widehat{\frac{x^4}{4} - \frac{x^2}{2}}^\mcn \quad x\in\er. 
\end{cases}
\eee
The analytic solution of \eqref{model_pr_ext} 
on $(0,1)\times (0, \frac{1}{2})$ is 
$$
u(t,x) = -\frac{14}{120} + \sum_{k=1}^{\infty} (-1)^{k+1}
\frac{12}{(k\pi)^4} e^{-\frac{t}{4}(k\pi)^{1.2}}\cos(k\pi x)
\quad x\in (0,1),\; t\in (0, 1),
$$
which has been computed in the grid points with a high accuracy 
to verify the convergence of the implicit Euler method. 
The results of the computations are summarized in Table \ref{res_imple}.
We computed the error $\|\ee_{\tau,h}\|_\infty$ of the approximation 
in maximum-norm for various time steps and discretization parameters
at the final time $t=1$. One can clearly see the first order convergence 
which was predicted by the theory, see Theorem \ref{main_thm}. 
The convergence rate was estimated in the consecutive refinement
steps using the formula 
$\log_2 \left(\frac{\|\ee_{2\tau,2h}\|_\infty}{\|\ee_{\tau,h}\|_\infty}\right)$.   
\begin{table}[ht]
\caption{Convergence results for the implicit Euler method 
in \eqref{impl_euler} applied to \eqref{model_pr_ext}}
\centering
    \begin{tabular}{| l | l | l | l |}
\hline\hline 
Grid parameter ($h$)  & Time step ($\tau$) & Error in $\|\cdot\|_\infty$-norm 
 & Convergence rate \\ [0.5ex] 
\hline 
1/2   & 1/2     & 1.4 $\cdot 10^{-2}$ & $\emptyset$ \\
1/4   & 1/4     & 1.5 $\cdot 10^{-2}$ & 0.1583 \\
1/8   & 1/8     & 1.2 $\cdot 10^{-2}$ & 0.3342 \\
1/16  & 1/16    & 7.8 $\cdot 10^{-3}$ & 0.6299 \\
1/32  & 1/32    & 4.5 $\cdot 10^{-3}$ & 0.8007 \\
1/64  & 1/64    & 2.4 $\cdot 10^{-3}$ & 0.8952 \\
1/128 & 1/128   & 1.3 $\cdot 10^{-3}$ & 0.9459 \\
1/256 & 1/256   & 6.3 $\cdot 10^{-4}$ & 0.9725 \\
1/512 & 1/512   & 3.2 $\cdot 10^{-4}$ & 0.9861 \\ [1ex] 
\hline
\end{tabular}
\label{res_imple}
\end{table}

The model of (fractional order) diffusion predicts that 
in case of homogeneous Neumann boundary conditions
the total mass should be preserved. Accordingly, in the 
numerical simulations above, the $l_1$ norm should be constant, 
which is an easy consequence of the fact, that the sum of 
the elements in the columns of $A_{\alpha,h,\infty}$ is zero. 

\begin{figure}[h!]
\centering
\includegraphics[scale=0.5]{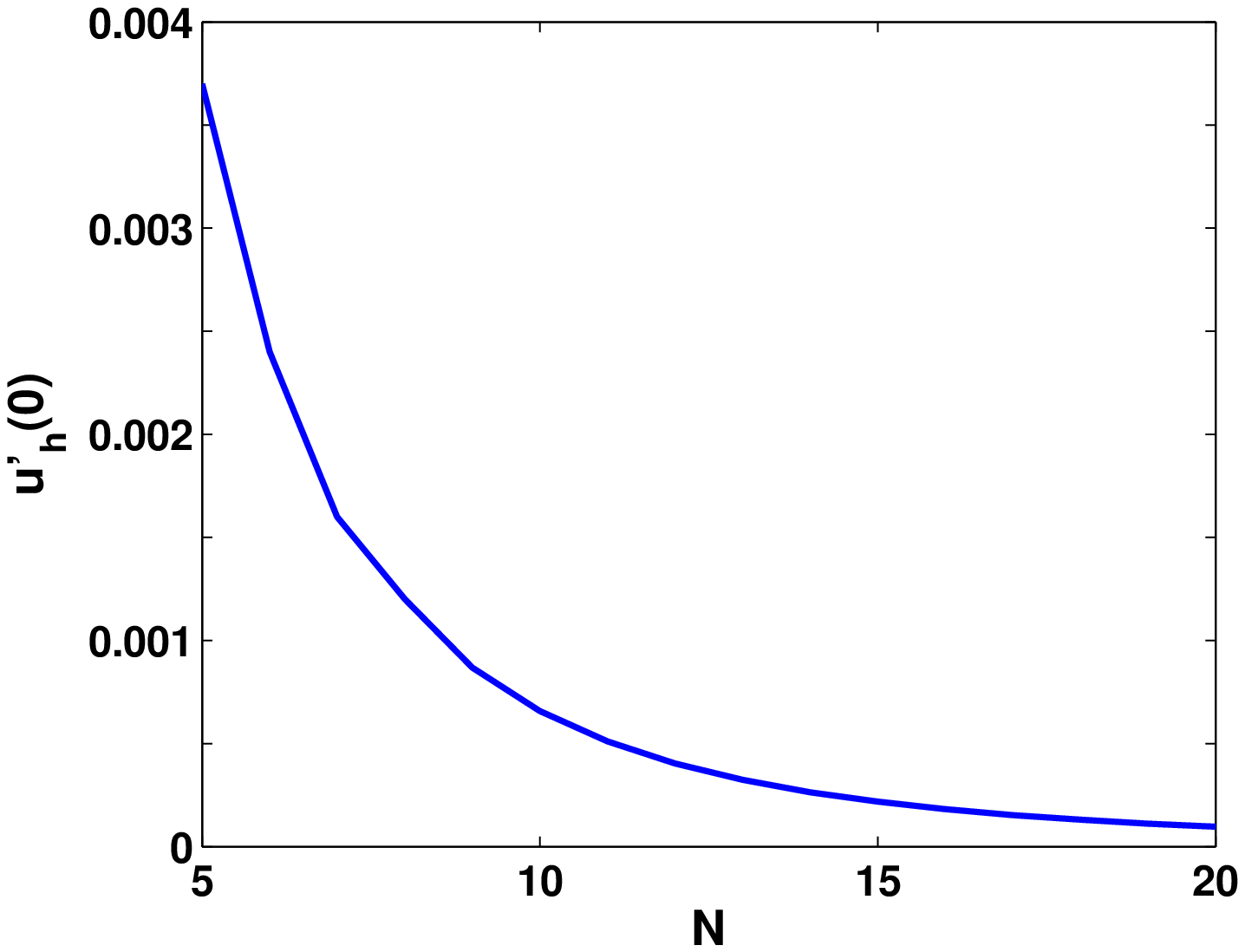}
\caption{The approximation \eqref{approx_n} of the derivative 
$\pa_x u(1,0)$ in the numerical simulations vs. the number of 
gridpoints.}
\label{neu_approx}
\end{figure}

Therefore, we examine the boundary conditions in course of the  
simulations. For this, we use the second order accurate approximation
\bee\label{approx_n}
\pa_x u(t,0) = \frac{1}{2h}(3u(t,0) - 4u(t,h) + u(t,2h)) 
\eee
and the results are shown in Figure \ref{neu_approx}. For the simplicity,
we applied the same number of grid points in the spatial and the time 
coordinates. This accurate approximation can be recognized as the 
numerical equivalent of Lemma \ref{compatible_lem}. 

\paragraph{An inhomogeneous model problem}
We second investigate the model problem
\bee\label{numex2_full}
\begin{cases}
\pa_t u(t,x) = 
0.5\cdot\pa_{|x|}^{1.6} u(t,x) + e^{-t}\cos \pi x
\quad x\in (0,1),\; t\in (0, 1) \\
u(0,x) = 2x^2 - \frac{4}{3} x^3\quad x\in (0,1)\\
\pa_x u (t,0) = \pa_x u (t,1) = 0\quad t\in (0, 1),
\end{cases}
\eee
which is converted to the well-posed extended problem
\bee\label{model_pr_ext2}
\begin{cases}
\pa_t u(t,x) = 
0.5\cdot\pa_{|x|}^{1.6} u(t,x) + e^{-t}\cos \pi x
\quad x\in (0,1),\; t\in (0, 1) \\
u(0,x) = \widehat{2x^2 - \frac{4}{3} x^3}^\mcn \quad x\in\er. 
\end{cases}
\eee
We split the original equation in 
\eqref{numex2_full} into two ones as 
\bees
\begin{cases}
\pa_t u_1(t,x) = 
0.5\cdot\pa_{|x|}^{1.6} u_1(t,x)
\quad x\in (0,1),\; t\in (0, 1) \\
u_1(0,x) = 2x^2 - \frac{4}{3} x^3\quad x\in (0,1)\\
\pa_x u_1 (t,0) = \pa_x u_1 (t,1) = 0\quad t\in (0, 1)
\end{cases}
\eees
and
\bees
\begin{cases}
\pa_t u_{2,t^*}(t,x) = 
0.5\cdot\pa_{|x|}^{1.6} u_{2,t^*}(t,x)\quad x\in (0,1),\; t\in (0, 1)\\
u_{2,t^*}(0,x) = e^{-t^*}\cos \pi x \quad x\in (0,1)\\
\pa_x u_{2,t^*} (t,0) = \pa_x u_{2,t^*} (t,1) = 0\quad t\in (0, 1),
\end{cases}
\eees
the analytic solution of which be given as 
$$
u_1(t,x) = 
\frac{1}{3} + \sum_{k=1}^{\infty}
 \frac{16}{(k\pi)^4} e^{-\frac{t}{2}(k\pi)^{1.6}}((-1)^k-1)\cos(k\pi x)
\quad x\in (0,1),\; t\in (0, 1),
$$
and
$$
u_{2,t^*} (t,x) = e^{-t^{1.6} - \frac{t\pi^{}}{2}}
\quad x\in (0,1),\; t\in (0, 1).
$$
These have been computed in the grid points with a high accuracy 
to verify the convergence of the implicit Euler method. \\

Then, according to the Duhamel's principle (see \cite{evans98}, p. 49)
we have 
$$
u(t,x) = u_1(t,x) + \int_0^t u_{2,t^*}(t-t^*,x)\:\mathrm{d}t^*.
$$
Accordingly, for the numerical solution at $t=1$ we compute first 
the approximation of $u_1(t,\cdot)$ at the grid points. Then with 
the same time steps and spatial accuracy we approximate 
$$
u_{2,0}(1,\cdot),
u_{2,\tau}(1-\tau,\cdot), u_{2,2\tau}(1-2\tau,x),\dots, u_{2,1}(0,x),
$$
where indeed, the last term is already given. Then a composite trapezoidal numerical integration
$$
u(1,x) \approx 
\frac{\tau}{2} u_{2,0}(1,\cdot)+
\tau (u_{2,\tau}(1-\tau,\cdot) + u_{2,2\tau}(1-2\tau,x)+ 
\dots + u_{2,1-\tau}(\tau,x)) +
\frac{\tau}{2} u_{2,1}(0,x)
$$
gives the desired approximation in $x$. 

The results of the computations are summarized in Table \ref{res_imple_2},
where in all cases $\tau = 1/1024$ such that the numerical integration does
not harm the predicted order of convergence. 

\begin{table}[ht]
\caption{Convergence results for the implicit Euler method 
in \eqref{impl_euler} applied to \eqref{numex2_full}}
\centering
    \begin{tabular}{| l | l | l | l |}
\hline\hline 
Grid parameter ($h$)  & Time step ($\tau$) & Error in $\|\cdot\|_\infty$-norm 
 & Convergence rate \\ [0.5ex] 
\hline 
1/2   & 1/2     & 6.2 $\cdot 10^{-2}$ & $\emptyset$ \\ 
1/4   & 1/4     & 3.7 $\cdot 10^{-2}$ & 0.7384 \\
1/8   & 1/8     & 2.1 $\cdot 10^{-2}$ & 0.8490 \\
1/16  & 1/16    & 1.1 $\cdot 10^{-3}$ & 0.9142 \\
1/32  & 1/32    & 5.7 $\cdot 10^{-3}$ & 0.9535 \\
1/64  & 1/64    & 2.9 $\cdot 10^{-3}$ & 0.9756 \\
1/128 & 1/128   & 1.5 $\cdot 10^{-3}$ & 0.9875 \\
1/256 & 1/256   & 7.3 $\cdot 10^{-4}$ & 0.9937 \\
1/512 & 1/512   & 3.7 $\cdot 10^{-4}$ & 0.9967 \\ [1ex] 
\hline
\end{tabular}
\label{res_imple_2}
\end{table}

\section*{Appendix}
\emph{Proof of Lemma \ref{lem_eig}:}
We first observe that the Neumann extension is the natural one for 
$\vv^k_h$ in the sense that 
\bee\label{natural_neu_ext}
\left[\widehat{\vv_h^k}^\mcn\right]_j = \cos kx_j \quad j\in\mathbb{Z}.
\eee
Then according to \eqref{A_def}, \eqref{natural_neu_ext}, 
\eqref{ezxj}, \eqref{trGLformula_left} and 
\eqref{trGLformula_right} we obtain 
\bee\label{long_eig}
\bali
&\left[A_{\alpha,h} \vv_h^k\right]_j
=
\sum_{l=0}^\infty \frac{C_\sigma}{h^\alpha} 
\left(
g_l \left[\widehat{\vv_h^k}^\mcn\right]_{j+l-1} 
+ g_l \left[\widehat{\vv_h^k}^\mcn\right]_{j-l+1}
\right)\\
&=
\frac{C_\sigma}{h^\alpha} \sum_{l=0}^\infty
g_l (\cos kx_{j+l-1} + \cos kx_{j-l+1})\\
&=
\frac{C_\sigma}{h^\alpha} \sum_{l=0}^\infty
g_l \left(\cos k\pi h \left(j+\frac{1}{2} + l-1)\right) 
     + \cos k\pi h\left(j+\frac{1}{2} - (l-1) \right) \right)\\
&=
\frac{2C_\sigma \cos kx_j}{h^\alpha} \sum_{l=0}^\infty
g_l \cos k\pi h(l-1) \\
&=
\frac{2C_\sigma \cos kx_j}{h^\alpha} \mathfrak{R}
\left(
\sum_{l=0}^\infty \exp\{-i k\pi h\} g_l \exp\{i k\pi hl\}
\right)\\
&=
\frac{2C_\sigma \cos kx_j}{h^\alpha} \mathfrak{R}
\left(
 \exp\{-i k\pi h\}\sum_{l=0}^\infty(-1)^l\binom{\alpha}{l}
 \exp\{i k\pi h l\}
\right)\\
&=
\frac{2C_\sigma \cos kx_j}{h^\alpha} \mathfrak{R}
\left(\exp\{-i k\pi h\} (1-\exp\{i k\pi h\})^\alpha\right)\\
&=
\frac{2C_\sigma \cos k x_j}{h^\alpha} \mathfrak{R}
\left( 
\exp\{-i k\pi h\}\cdot 
2^\alpha\sin^\alpha \frac{k\pi h}{2}\cdot
\left(\cos\frac{\alpha}{2}(\pi- k\pi h) - 
  i \sin\frac{\alpha}{2}(\pi- k\pi h)\right)
\right)\\
&=
\frac{2C_\sigma \cos kx_j}{h^\alpha}
2^\alpha\sin^\alpha \frac{k\pi h}{2}\cdot 
\left( \cos k\pi h \cos\frac{\alpha}{2}(\pi-k\pi h) 
-\sin k\pi h \sin\frac{\alpha}{2}(\pi-k\pi h)
\right)\\
&=
\frac{2C_\sigma \cos kx_j}{h^\alpha}
2^\alpha\sin^\alpha \frac{k\pi h}{2}\cdot 
\cos \left(k\pi h + \frac{\alpha}{2}(\pi-k\pi h)\right)
\eali
\eee
where the have used the identity 
\[
\bali
&(1-\exp\{i k\pi h\})^\alpha =
\left(2\sin\frac{k\pi h}{2}\cdot
\left(\sin\frac{k\pi h}{2} - i \cos\frac{k\pi h}{2}\right)\right)^\alpha\\
&=
2^\alpha\sin^\alpha \frac{k}{2}\pi h\cdot
\left(\cos\left(\frac{\pi}{2}-\frac{k\pi h}{2}\right) - 
  i \sin\left(\frac{\pi}{2}-\frac{k\pi h}{2}\right)\right)^\alpha\\
&=
2^\alpha\sin^\alpha \frac{k\pi h}{2}\cdot
\left(\cos\frac{\alpha}{2}(\pi-k\pi h) - 
  i \sin\frac{\alpha}{2}(\pi-k\pi h) \right).
\eali
\]
The definition of $C_\sigma$ gives then the statement in the lemma.
\quad$\square$

\emph{Proof of Lemma \ref{eig_cos_lem}:}
In the proof, we use the identities 
\bee\label{grads}
\bali
\int_0^\infty x^{n-1} \cos bx \dx &= \frac{\Gamma(n)}{b^n} \cos\frac{n\pi}{2},\\
\int_0^\infty x^{n-1} \sin bx \dx &= \frac{\Gamma(n)}{b^n} \sin\frac{n\pi}{2},
\eali
\eee
which can be found in \cite{gradstein00}, 3.761/9.

Observe that the even extension of the $\cos(k\pi\cdot)|_{(0,1)}$ function
to the real axis is the $\cos(k\pi\cdot)$ function itself. On the other 
hand, as it was pointed out in \cite{yang10}, we can differentiate the 
integrals in the Riemann--Liouville formula to obtain 
\bee\label{last_eq}
\pa^\alpha_{|x|} \cos(k\pi x) =
-\frac{\sigma\cdot(k\pi)^2}{2\cos\left(\frac{\alpha\pi}{2}\right)}
\left(
-\int_\infty^x\frac{\cos k\pi s}{(x-s)^{\alpha -1}}\ds 
-\int_\infty^x\frac{\cos k\pi s}{(s-x)^{\alpha -1}}\ds
\right).
\eee 
Using \eqref{grads}, the first term can be rewritten as
$$
\bali
&\int_{-\infty}^x\frac{\cos k\pi s}{(x-s)^{\alpha -1}}\ds
=
\int_0^\infty \frac{\cos k\pi (x-y)}{y^{\alpha -1}}\dy\\
&=
\cos k\pi x \int_0^\infty  \frac{\cos k\pi y}{y^{\alpha -1}}\dy
+
\sin k\pi x \int_0^\infty  \frac{\sin k\pi y}{y^{\alpha -1}}\dy\\
&=
\frac{\Gamma(2-\alpha)}{(k\pi)^{2-\alpha}}
\left(\cos k\pi x\cos\frac{\pi(2-\alpha)}{2} + 
\sin k\pi x\sin\frac{\pi(2-\alpha)}{2}\right).
\eali
$$
A similar computations gives that 
$$
\bali
&\int^{\infty}_x\frac{\cos k\pi s}{(s-x)^{\alpha -1}}\ds
=
\int_0^\infty \frac{\cos k\pi (x+y)}{y^{\alpha -1}}\dy\\
&=
\frac{\Gamma(2-\alpha)}{(k\pi)^{2-\alpha}}
\left(\cos k\pi x\cos\frac{\pi(2-\alpha)}{2} - 
\sin k\pi x\sin\frac{\pi(2-\alpha)}{2}\right).
\eali
$$
Therefore, the equality in \eqref{last_eq} can be rewritten as 
\bees
\bali
\pa^\alpha_{|x|} \cos(k\pi x) &=
\frac{\sigma\cdot(k\pi)^2}
{2\cos\left(\frac{\alpha\pi}{2}\right)\Gamma(2-\alpha)}
\frac{\Gamma(2-\alpha)}{(k\pi)^{2-\alpha}}
2 \cos k\pi x\cos\frac{\pi(2-\alpha)}{2}\\
&=
\frac{\sigma\cdot(k\pi)^\alpha}{2\cos\left(\frac{\alpha\pi}{2}\right)}
2 \cos k\pi x\cos\frac{\pi(2-\alpha)}{2}
=
-\sigma\cdot(k\pi)^\alpha \cos k\pi x.
\eali
\eees
On the other hand, the system $\left\{ \cos k\pi x\right\}_{k=1}^{\infty}$
is complete in $L_2 (0,1)$ such that no further eigenfunctions can exist. 
\quad $\square$

\section*{Acknowledgements}
The authors acknowledge the financial support of the 
Hungarian National Research Fund OTKA (grants K104666 and K112157).

\bibliography{ferenc_bib.bib}

\def\cprime{$'$}
\begin{thebibliography}{10}

\bibitem{adams03}
R.~A. Adams and J.~J. Fournier.
\newblock {\em Sobolev spaces}.
\newblock Academic Press, Amsterdam, second edition, 2003.
\newblock Pure and Applied Mathematics, Vol. 140.

\bibitem{benson00}
D.~A. Benson, S.~W. Wheatcraft, and M.~M. Meerschaert.
\newblock Application of a fractional advection-dispersion equation.
\newblock {\em Water Resources Research}, 36(6):1403--1412 (electronic), 2000.

\bibitem{bonito13}
A.~Bonito and J.~Pasciak.
\newblock Numerical approximation of fractional powers of elliptic operators.
\newblock Submitted, arXiv:1307.0888.

\bibitem{deng14}
W.~H. Deng and M.~Chen.
\newblock Efficient numerical algorithms for three-dimensional fractional
  partial diffusion equations.
\newblock {\em J. Comp. Math.}, 32(4):371--391, 2014.

\bibitem{du12}
Q.~Du, M.~Gunzburger, R.~Lehoucq, and K.~Zhou.
\newblock Analysis and approximation of nonlocal diffusion problems with volume
  constraints.
\newblock {\em SIAM Rev.}, 54(4):667--696, 2012.

\bibitem{du13}
Q.~Du, M.~Gunzburger, R.~B. Lehoucq, and K.~Zhou.
\newblock A nonlocal vector calculus, nonlocal volume-constrained problems, and
  nonlocal balance laws.
\newblock {\em Math. Mod. Meth. Appl. Sci.}, 23(03):493--540, 2013.

\bibitem{edwards07}
A.~M. Edwards, R.~A. Phillips, N.~W. Watkins, M.~P. Freeman, E.~J. Murphy,
  V.~Afanasyev, S.~V. Buldyrev, M.~G.~E. da~Luz, E.~P. Raposo, H.~E. Stanley,
  and G.~M. Viswanathan.
\newblock Revisiting {L}\'evy flight search patterns of wandering albatrosses,
  bumblebees and deer.
\newblock {\em Nature}, 449:1044--1048, 2007.

\bibitem{evans98}
L.~C. Evans.
\newblock {\em Partial differential equations}, volume~19 of {\em Graduate
  Studies in Mathematics}.
\newblock American Mathematical Society, Providence, RI, 1998.

\bibitem{gradstein00}
I.~S. Gradshteyn and I.~M. Ryzhik.
\newblock {\em Table of integrals, series, and products}.
\newblock Academic Press Inc., San Diego, CA, sixth edition, 2000.
\newblock Translated from the Russian, Translation edited and with a preface by
  Alan Jeffrey and Daniel Zwillinger.

\bibitem{hilfer08}
R.~Hilfer.
\newblock Threefold introduction to fractional derivatives.
\newblock In R.~Klages, G.~Radons, and I.~Sokolov, editors, {\em Anomalous
  Transport: Foundations and Applications}, pages 17--73. Wiley-VCH Verlag GmbH
  $\&$ Co. KGaA, Weinheim, 2008.

\bibitem{huang13}
J.-F. Huang, N.~Nie, and Y.-F. Tang.
\newblock A second order finite difference-spectral method for space fractional
  diffusion equations.
\newblock {\em Sci. Chin. Math.}, 57(6):1303--1317, 2014.

\bibitem{ilic05}
M.~Ilic, F.~Liu, I.~Turner, and V.~Anh.
\newblock Numerical approximation of a fractional-in-space diffusion equation.
  {I}.
\newblock {\em Fract. Calc. Appl. Anal.}, 8(3):323--341, 2005.

\bibitem{kilbas06}
A.~A. Kilbas, H.~M. Srivastava, and J.~J. Trujillo.
\newblock {\em Theory and applications of fractional differential equations},
  volume 204 of {\em North-Holland Mathematics Studies}.
\newblock Elsevier Science B.V., Amsterdam, 2006.

\bibitem{meerschaert04}
M.~M. Meerschaert and C.~Tadjeran.
\newblock Finite difference approximations for fractional advection-dispersion
  flow equations.
\newblock {\em J. Comput. Appl. Math.}, 172(1):65--77, 2004.

\bibitem{miller93}
K.~S. Miller and B.~Ross.
\newblock {\em An introduction to the fractional calculus and fractional
  differential equations}.
\newblock A Wiley-Interscience Publication. John Wiley \& Sons Inc., New York,
  1993.

\bibitem{nochetto14}
R.~Nochetto, E.~Otárola, and A.~Salgado.
\newblock A pde approach to fractional diffusion in general domains: A priori
  error analysis.
\newblock {\em Foundations of Computational Mathematics}, 2014.

\bibitem{podlubny99}
I.~Podlubny.
\newblock {\em Fractional differential equations}, volume 198 of {\em
  Mathematics in Science and Engineering}.
\newblock Academic Press Inc., San Diego, CA, 1999.
\newblock An introduction to fractional derivatives, fractional differential
  equations, to methods of their solution and some of their applications.

\bibitem{samko93}
S.~G. Samko, A.~A. Kilbas, and O.~I. Marichev.
\newblock {\em Fractional integrals and derivatives}.
\newblock Gordon and Breach Science Publishers, Yverdon, 1993.
\newblock Theory and applications, Edited and with a foreword by S. M.
  Nikol{\cprime}ski{\u\i}, Translated from the 1987 Russian original, Revised
  by the authors.

\bibitem{shen08}
S.~Shen, F.~Liu, V.~Anh, and I.~Turner.
\newblock The fundamental solution and numerical solution of the {R}iesz
  fractional advection-dispersion equation.
\newblock {\em IMA J. Appl. Math.}, 73(6):850--872, 2008.

\bibitem{szymczak03}
P.~Szymczak and A.~J.~C. Ladd.
\newblock Boundary conditions for stochastic solutions of the
  convection-diffusion equation.
\newblock {\em Phys. Rev. E}, 68:036704, Sep 2003.

\bibitem{tadjeran07}
C.~Tadjeran and M.~M. Meerschaert.
\newblock A second-order accurate numerical method for the two-dimensional
  fractional diffusion equation.
\newblock {\em J. Comput. Phys.}, 220(2):813--823, 2007.

\bibitem{tadjeran06}
C.~Tadjeran, M.~M. Meerschaert, and H.-P. Scheffler.
\newblock A second-order accurate numerical approximation for the fractional
  diffusion equation.
\newblock {\em J. Comput. Phys.}, 213(1):205--213, 2006.

\bibitem{tian12}
W.~Tian, H.~Zhou, and W.~Deng.
\newblock A class of second order difference approximation for solving space
  fractional diffusion equations.
\newblock {\em Math. Comp.}, 2014.
\newblock to appear, {\tt arXiv:1201.5949[math.NA]}.

\bibitem{treumann97}
R.~A. Treumann.
\newblock Theory of super-diffusion for the magnetopause.
\newblock {\em Geophys. Res. Lett.}, 24(14):1727--1730, 1997.

\bibitem{wang12}
H.~Wang and T.~S. Basu.
\newblock A fast finite difference method for two-dimensional space-fractional
  diffusion equations.
\newblock {\em SIAM J. Sci. Comput.}, 34(5):A2444--A2458, 2012.

\bibitem{yang10}
Q.~Yang, F.~Liu, and I.~Turner.
\newblock Numerical methods for fractional partial differential equations with
  {R}iesz space fractional derivatives.
\newblock {\em Appl. Math. Model.}, 34(1):200--218, 2010.

\bibitem{zhou13}
H.~Zhou, W.~Tian, and W.~Deng.
\newblock Quasi-compact finite difference schemes for space fractional
  diffusion equations.
\newblock {\em J. Sci. Comput.}, 56(1):45--66, 2013.

\end{thebibliography}
\bibliographystyle{abbrv}

\end{document}